\documentclass[12pt]{amsart}

\usepackage{amsfonts,amssymb,amsmath}
\usepackage{graphicx}
\usepackage{psfrag}
\usepackage{color,xcolor}
\usepackage{verbatim}
\usepackage{epstopdf}
\usepackage[text={6in,9in},centering]{geometry}
\usepackage[colorlinks,citecolor=red,pagebackref,hypertexnames=false]{hyperref}
\newtheorem{thm}{Theorem}[section]
\newtheorem{lem}[thm]{Lemma}

\theoremstyle{remark}
\newtheorem{rem}[thm]{Remark}
\newtheorem{conj}[thm]{Conjecture}

\numberwithin{equation}{section}

\newcommand{\R}{{\mathbb R}}

\newcommand{\Z}{{\mathbb Z}}

\newcommand{\N}{{\mathcal N}}

\newcommand{\D}{{\mathcal D}}

\newcommand{\T}{{\mathcal T}}
\newcommand{\K}{{\mathcal K}}

\def\P{{\mathcal P}}

\def\bfn{{\mathbf n}}
\def\bfm{{\mathbf m}}
\def\bfa{{\mathbf a}}
\def\bfb{{\mathbf b}}

\def\id{{\operatorname{id}}}

\def\bfx{{\mathbf x}}

\let\rho\varrho

\newcommand{\bfeta}{{\boldsymbol{\eta}}}
\newcommand{\bftau}{{\boldsymbol{\tau}}}
\newcommand{\bfsig}{{\boldsymbol{\sigma}}}

\def\updimB{{\overline{\dim} _{\rm M}}}
\def\lowdimB{{\underline{\dim} _{\rm M}}}
\def\dimB{\dim_{\rm M}}

\def\wdt{\widetilde}
\def\vep{\varepsilon}
\def\leb{\mathcal{L}}
\newcommand{\card}{{\rm card\,}}

\newcommand{\lip}{{\rm Lip\,}}

\begin{document}

\title{Minkowski dimension of iterated sum-product sets}

\author{Yifei Gu}
\address{Department of Mathematics, East China Normal University,  Shanghai 200241, P. R. China}
\email{52275500012@stu.ecnu.edu.cn}

\author{Lai Jiang}
\address{School of Science, Shanghai Maritime University, Shanghai 201306, P. R. China}
\email{jianglai@shmtu.edu.cn}


\thanks{Corresponding author: Lai Jiang}


\date{}


\begin{abstract}
For $m\in\mathbb{Z}^+$ and $\tau>0$, let $F_{m,\tau}$ denote the $m$-fold arithmetic sum of
$$
\{n^{-\tau}:n\in\mathbb{Z}^+\}\cup\{0\}.
$$
We determine the Minkowski dimensions of finite products
$$
\mathcal{P}_{\mathbf{m},\boldsymbol{\tau}}
:=\prod_{i=1}^s F_{m_i,\tau_i},
$$
where both the number of summands and the decay exponent may vary from factor to factor.
We also study the Taylor sets
$$
\mathcal{T}_{k,m}:=\sum_{i=1}^{k}\prod_{j=1}^{i}F_{m,1}.
$$
For $k = 2$, we determine the Minkowski dimension exactly and compare the Taylor set with a product set of the same dimension.
For $k \geq 3$, we give a lower bound for the Minkowski dimension.
\end{abstract}

\maketitle

\section{Introduction}

\subsection{Motivation and background}

For subsets $A,B\subset\R$, write
$$
A+B=\{a+b:a\in A,\ b\in B\},
\qquad
AB=\{ab:a\in A,\ b\in B\}.
$$
For a finite set $A$ of integers, the Erd\H{o}s--Szemer\'edi problem asks how large at least one of $A+A$ and $AA$ must be in terms of the size of $A$ \cite{ES83}.
For finite sets of real numbers, see \cite{E97,S09} for quantitative bounds.
The product of the middle-third Cantor set with itself contains an interval \cite{XJZP19}; for the topology and Lebesgue measure of products of uniform Cantor sets, see \cite{GJXZ21}.
For arithmetic sums of self-similar sets on the line, the Hausdorff dimension can depend on arithmetic relations among the contraction ratios \cite{J16,PS09}.

Iterating addition raises a further question: how does the dimension of a sumset grow with the number of summands? For compact subsets of the line, Schmeling and Shmerkin \cite{SS10} showed that the Hausdorff dimensions of successive sumsets can follow any nondecreasing sequence in $[0,1]$.
Fraser, Howroyd and Yu \cite{FHY19} showed that these dimensions tend to one when the original set is compact and uniformly perfect.

We study the Minkowski dimensions of a class of iterated sum-product sets. Let
$$
F=\{n^{-1}:n\in\Z^+\}\cup\{0\},
$$
and let $F_m$ denote the $m$-fold arithmetic sum of $F$. 
Every set obtained from $F$ by finitely many sums and products is countable, and therefore has Hausdorff dimension zero. Minkowski dimension can distinguish such sets through the growth of their covering numbers \cite{F90}.
Sums of unit fractions also arise in the theory of number expansions \cite{DK02}. 
The growth of digits has been studied for continued fractions, Engel expansions and Sylvester series \cite{FT21,KW04,LR16,LW01}.
Fraser, Howroyd and Yu \cite{FHY19} showed that the lower Minkowski dimension of $F_m$ is at least $1-2^{-m}$.
In their study of Egyptian fraction expansions, Chen, Jiang and Wu \cite{CJW} proved that
\begin{equation}\label{eq:fm-dim}
\dimB(F_m)=1-2^{-m}.
\end{equation}
Thus each additional copy of $F$ halves the gap between the dimension and one. The gap remains positive for every finite $m$.

We wonder how multiplication affects the Minkowski dimension of these sumsets.
We first determine the dimension of finite products of
$F_{m_1},\ldots,F_{m_s}$, allowing more generally for different decay exponents in the factors. 
We then study sums of such products and compare the dimension of $F_m+F_pF_q$ with that of $F_{m+p}F_q$ when $p \geq q$.

\subsection{Notation and basic definitions}

For $\tau>0$, set
$$
F^{(\tau)}=\{n^{-\tau}:n\in\Z^+\}\cup\{0\}, \qquad F_{0,\tau}=\{0\},
$$
and let $F_{m,\tau}$ be the $m$-fold arithmetic sum of $F^{(\tau)}$ for
$m\ge1$. In particular, $F^{(1)}=F$ and $F_{m,1}=F_m$.

For $s\in\Z^+$, set $\Sigma^s=(\Z^+)^s$. 
Given $\bfm=(m_1,\ldots,m_s)\in\Sigma^s$ and $\bftau=(\tau_1,\ldots,\tau_s)\in(0,\infty)^s$, define the generalized product set
$$
\P_{\bfm,\bftau} =\prod_{i=1}^{s}F_{m_i,\tau_i}.
$$
When all decay exponents are one, we write $\P_{\bfm}=F_{m_1}\cdots F_{m_s}$.

For $k,m\in\Z^+$, define the Taylor set
$$
\T_{k,m}
=\sum_{i=1}^{k}\prod_{j=1}^{i}F_m
=\bigg\{\sum_{i=1}^{k}\prod_{j=1}^{i}x_{i,j}:x_{i,j}\in F_m \bigg\},
$$
with each $x_{i,j}$ chosen independently. 
In particular, $\T_{1,m}=F_m$ and $\T_{2,m}=F_m+F_mF_m$.

For a nonempty bounded set $\Omega\subset\R$, let $\N_r(\Omega)$ be the least number of closed intervals of length
$r$ needed to cover $\Omega$. 
Its upper and lower Minkowski dimensions are
\begin{equation}\label{eq:min-dim-cover}
\updimB(\Omega)=\limsup_{r\to0^+}\frac{\log\N_r(\Omega)}{-\log r}, \qquad
\lowdimB(\Omega)=\liminf_{r\to0^+}\frac{\log\N_r(\Omega)}{-\log r}.
\end{equation}
When the two values coincide, we denote their common value by
$\dimB(\Omega)$. Writing $\leb$ for Lebesgue measure on $\R$ and
$\Omega_r=\{x\in\R:|x-y|\leq r\text{ for some }y\in\Omega\}$, we have the
equivalent formulas
\begin{equation}\label{eq:min-dim-leb}
\updimB(\Omega) =1-\liminf_{r\to0^+}\frac{\log\leb(\Omega_r)}{\log r},\qquad
\lowdimB(\Omega) =1-\limsup_{r\to0^+}\frac{\log\leb(\Omega_r)}{\log r}.
\end{equation}
See \cite{F90} for standard properties of Minkowski dimension.

\subsection{Main results}

Two parameters determine the dimension of a finite product of unit-fraction sumsets.

\begin{thm}\label{thm:prod}
Let $s\in\Z^+$ and $\bfm\in\Sigma^s$. 
With $M=\max_i m_i$ and $K=\sum_{i=1}^s(m_i-1)$, we have
$$
\dimB(\P_{\bfm}) = 1-\frac{1}{1+(2^M-1)2^{K-M+1}}.
$$
\end{thm}

For two factors, this gives
$$
\dimB(F_pF_q) = 1-\frac{1}{1+(2^{\max\{p,q\}}-1)2^{\min\{p,q\}-1}}.
$$

To allow different decay exponents in the factors, set $\kappa_\tau=(1+\tau)/\tau$ for $\tau>0$. 
For a single summand, $\dimB(F_{1,\tau})=\frac{1}{1+\tau}$ by the argument of \cite[Example~3.5]{F90}. 
See also \cite{MZ99,ZM05} for general formulas for decreasing convex sequences.
For the $m$-fold sum, Chen, Jiang and Wu \cite[Remark~1.5]{CJW} give
$$
\dimB(F_{m,\tau})=1-\kappa_\tau^{-m}.
$$
For products, define
\begin{equation}\label{eq:generalized-D-short}
D_0(\bfm,\bftau) =\prod_{j=1}^s\kappa_{\tau_j}^{m_j-1},\qquad
D_i(\bfm,\bftau) =1+\big(\kappa_{\tau_i}^{m_i}-1\big)  \prod_{j\ne i}\kappa_{\tau_j}^{m_j-1}, \quad 1\leq i\leq s,
\end{equation}
and set
\begin{equation}\label{eq:generalized-D}
D(\bfm,\bftau)=\max_{0\leq i\leq s}D_i(\bfm,\bftau).
\end{equation}

\begin{thm}\label{thm:gen}
Let $s\in\Z^+$, $\bfm\in\Sigma^s$ and
$\bftau\in(0,\infty)^s$. Then
$$
\dimB(\P_{\bfm,\bftau}) = 1-\frac{1}{D(\bfm,\bftau)}.
$$
\end{thm}

The quadratic Taylor set $\T_{2,m}=F_m+F_mF_m$ is the simplest case of adding a sumset to a product.

\begin{thm}\label{thm:true-taylor}
For every $m\in\Z^+$,
$$
\dimB(\T_{2,m}) =1-\frac{1}{1+2^{3m-1}-2^{m-1}}.
$$
\end{thm}

Combining Theorems~\ref{thm:prod} and~\ref{thm:true-taylor} gives $\dimB(\T_{2,m})=\dimB(F_{2m}F_m)$. 
For $m,p,q\in\Z^+$ with $p\geq q$, the proof of Theorem~\ref{thm:true-taylor} extends to $F_m+F_pF_q$ and gives
\begin{equation}\label{eq:dim-equal-prod-taylor}
\dimB(F_m+F_pF_q)
=1-\frac{1}{1+(2^{m+p}-1)2^{q-1}}
=\dimB(F_{m+p}F_q).
\end{equation}
Thus adding $F_m$ has the same effect on dimension as increasing the number of summands in the larger factor by $m$. 
We next compare the two sets under inclusion.

\begin{thm}\label{thm:non-equal}
For every $m,p,q\in\Z^+$ with $p,q\ge2$, the sets $F_m+F_pF_q$ and $F_{m+p}F_q$ are incomparable under inclusion.
\end{thm}

The Taylor sets of higher degree lead to the following conjecture.

\begin{conj}\label{conj:true-taylor-all}
For every $k,m\in\Z^+$, the Minkowski dimension of $\T_{k,m}$ exists and
$$
\dimB(\T_{k,m})=1-\frac{1}{d_{k,m}},
$$
where
\begin{equation}\label{eq:true-taylor-D}
d_{k,m}
=2^{km+k(k-1)(m-1)/2}
-\sum_{i=1}^{k-1} 2^{(k-i)(k-i-1)(m-1)/2+i-1}
  \big(2^{(k-i)(m-1)}-1\big).
\end{equation}
\end{conj}

Formula \eqref{eq:fm-dim} and Theorem~\ref{thm:true-taylor} establish the cases $k=1$ and $k=2$, respectively.
The conjecture also holds for $m=1$, since $FF=F$, $\T_{k,1}=F_k$ and $d_{k,1}=2^k$.
For $k\geq3$ and $m\geq2$, we prove the lower bound
$$
\lowdimB(\T_{k,m}) \geq 1-\frac{1}{d_{k,m}}.
$$

\subsection{Proof strategy and organization}

The proof of Theorem~\ref{thm:gen} uses two types of approximation for the lower bound. 
Approximation near $1$ yields $D_0$, while starting near $0$ with the $i$-th factor yields $D_i$. 
The maximum of these exponents gives the lower bound. 
For the upper bound, we decompose each factor according to its largest summand and estimate the covering numbers of the products of the resulting pieces. 
This keeps track of the dependence between terms introduced by multiplication.

For Taylor sets, we obtain the conjectured lower bound by constructing
successively finer approximations through induction on the degree. 
For the upper bound in degree two, we group the representations of elements of $F_m$ by their initial denominators in nondecreasing order.
Within each group, the corresponding unit fractions have a fixed sum, so it suffices to cover $F_mF_m$ plus the sum of the remaining unit fractions.

Section~2 proves Theorem~\ref{thm:gen}, deduces Theorem~\ref{thm:prod} and illustrates the product formula with $F_2F_2$. 
Section~3 establishes the exact formula for quadratic Taylor sets and the lower bound in higher degrees. 
Section~4 proves the non-containment result.

\section{Products of generalized sumsets}\label{sec:generalized-products}
In this section, we study generalized product sets arising from finite products of generalized sumsets and determine their Minkowski dimensions.
Recall that for every $\tau>0$, $\kappa_\tau := (1+\tau)/\tau$.

\subsection{The lower bound}

\begin{lem}\label{lem:gen-approx}
Fix an integer $m\geq 0$ and $\tau>0$.
There exists a constant $C>0$ depending only on $m$ and $\tau$ such that
for every $x\in(0,1)$, there exists $y\in F_{m, \tau}$ satisfying
$$
0<x-y\leq C x^{\kappa_\tau^m}.
$$
\end{lem}

\begin{proof}
Write $A_\tau := \tau 2^{\tau+1}$.
For any $z\in(0,1)$, there exists $n\geq 2$ such that
$
n^{-\tau}<z\leq(n-1)^{-\tau}.
$
By the mean value theorem, we have
\begin{equation}\label{eq:tau-single-low-est}
0<z-n^{-\tau}
\leq (n-1)^{-\tau}-n^{-\tau}
\leq \tau(n-1)^{-\tau-1}
\leq \tau 2^{\tau+1}n^{-\tau-1}
< A_\tau z^{\kappa_\tau}.
\end{equation}

We prove the lemma by induction.
For $m=0$, the lemma holds by taking $C=1$ and $y=0$.
Now we suppose $m\geq 1$.
Starting with $z_0=x$, we apply \eqref{eq:tau-single-low-est} successively.
For $1\leq j\leq m$, choose
$n_j\in\Z^+$ and put $z_j=z_{j-1}-n_j^{-\tau}$ so that
$0<z_j\leq A_\tau z_{j-1}^{\kappa_\tau}$. It follows that
$$
0<z_m
\leq A_\tau z_{m-1}^{\kappa_\tau}
\leq A_\tau (  A_\tau z_{m-2}^{\kappa_\tau})^{\kappa_\tau}
\leq \cdots
\leq A_\tau^{1+\kappa_\tau+\cdots+\kappa_\tau^{m-1}} x^{\kappa_\tau^m}.
$$
Let $C :=  A_\tau^{1+\kappa_\tau+\cdots+\kappa_\tau^{m-1}} $ and
$$
y=\sum_{j=1}^m n_j^{-\tau}  \in F_{m , \tau}
$$
Then we have $0 < x-y=z_m \leq C x^{\kappa_\tau^m}$.
Therefore the conclusion follows.
\end{proof}

Lemma~\ref{lem:gen-approx} leads to the following multiplicative refinement.

\begin{lem}\label{lem:gen-gap-pass}
Fix an integer $m\in \Z^+$ and $\tau , \gamma, C>0$.
There exists $C_{m,\tau}>0$ such that for every $x,y \in (0,1)$ with $0<x-y\leq C x^{1+\gamma}$,
there exists $u \in  F_{m, \tau}$ satisfying
$$
0\leq x-yu \leq (2C)^{\kappa_\tau^{m-1}} C_{m,\tau} x^{1+\gamma (\kappa_\tau)^{m-1}}.
$$

\end{lem}

\begin{proof}
For $m=1$, we take $C_{1,\tau}=1$ and $u=1$. Then by the assumption,
$$
0<x-yu = x-y\leq C x^{1+\gamma}
\leq 2 C x^{1+\gamma} .
$$
Since $m=1$ and $C_{1,\tau}=1$, this is the desired estimate.

Now we suppose $m\geq 2$.
Write $\vep=(x-y)/y$ and $\eta = \kappa_\tau^{m-1}$.
If $x < 2y$, then $0<\vep < 1$ and
$$
\vep=\frac{x-y}{y}
\leq \frac{2(x-y)}{x}
\leq 2C x^\gamma.
$$
By Lemma~\ref{lem:gen-approx}, there exists
$z\in F_{m-1,\tau}$ such that $0<\vep-z \leq \wdt{C}_{m-1,\tau}\vep^\eta$.
Since $u=1+z\in F_{m, \tau}$, we obtain
$$
0<x-yu
=y(\vep-z)
\leq \wdt{C}_{m-1,\tau}x\vep^\eta
\leq \wdt{C}_{m-1,\tau}(2C )^\eta
x^{1+\gamma \eta}.
$$

If $x \geq 2y$, take $u=1$. 
Notice that $ x \leq 2(  x-y) \leq 2 C x^{1+\gamma}$, which implies that $  2C x^{\gamma} \geq 1$.
Therefore
$$
0 < x-yu \leq x \leq (2C )^\eta x^{1+\gamma \eta}.
$$
Thus the conclusion holds by taking $C_{m, \tau }=  \max\{1, \wdt{C}_{m-1,\tau}\}$.
\end{proof}

\begin{lem}\label{lem:gen-gap-one}
Let $s\geq 2$, $\bfm \in \Sigma^s$ and $\bftau \in (0,\infty)^s$.
Then there exists $C >0$ such that for every $x \in (0,1)$, there exists $t \in \P_{\bfm, \bftau}$ such that
$$
	0 \leq x - t \leq C x^{D_1(\bfm, \bftau)}
$$

\end{lem}

\begin{proof}
Without loss of generality, we assume that $x \notin \P_{\bfm, \bftau}$.
Set $\gamma_1 := (\kappa_{\tau_1})^{m_1}-1$.
Define $\gamma_j := \gamma_{j-1} (\kappa_{\tau_j})^{m_j-1}$ for each $2 \leq j \leq s$.
By definition, $\gamma_s = D_1(\bfm, \bftau) -1 $.
Now we prove by induction that for each $1 \leq j \leq s$, there exists $ t_j \in F_{m_1, \tau_1} \cdots F_{m_j, \tau_j} $ such that
\begin{equation}\label{eq:gene-lower-gap-single-tmp}
	0 < x - t_j \leq C_j x^{\gamma_j +1}.
\end{equation}
Here $C_j>0$ is a constant independent of $ x $.
Let $C_1>0$ be the constant mentioned in Lemma~\ref{lem:gen-approx}.
Thus there exists $ t_1\in F_{m_1 , \tau_1}$ such that
$$
	0 < x- t_1\leq C_1 x^{ (\kappa_{\tau_1} )^{m_1}} = C_1 x^{\gamma_1+1} .
$$
This implies that \eqref{eq:gene-lower-gap-single-tmp} holds for $j=1$.

Now we assume that \eqref{eq:gene-lower-gap-single-tmp} holds for all $k<j$.
By Lemma~\ref{lem:gen-gap-pass}, there exists $u_j\in   F_{m_j, \tau_j}$ such that
$$
0 \leq x - t_{j-1} u_j \leq  C_j x^{ \gamma_{j-1} (\kappa_{\tau_j})^{m_j-1} +1} \leq
C_j x^{\gamma_j +1}.
$$
Here $C_j$ is a constant depending only on the fixed parameters.
Let $ t_j := t_{j-1} u_j $.
Then we have $ t_j \in F_{m_1, \tau_1} \cdots F_{m_j, \tau_j} \subset \P_{\bfm, \bftau} $.
Therefore $ x- t_j>0 $ follows from $x \notin \P_{\bfm, \bftau}$ and \eqref{eq:gene-lower-gap-single-tmp} holds for $j$.
Thus the conclusion follows by taking $t=t_s$ and $C=C_s$.
\end{proof}

The argument of the next lemma is similar to that in Lemma~\ref{lem:gen-gap-one}.
Here we focus on points near $1$.

\begin{lem}\label{lem:gen-gap-zero}
Let $s\geq 2$, $\bfm \in \Sigma^s$ and $\bftau \in (0,\infty)^s$. 
Then there exists $C>0$ such that for every $x\in(0,1)$, there exists $t\in\P_{\bfm,\bftau}$ such that
$$
0\leq 1+x-t\leq Cx^{D_0(\bfm,\bftau)}.
$$
\end{lem}

\begin{proof}
For each $1\leq j\leq s$, we set $\beta_j:=(\kappa_{\tau_j})^{m_j-1}$ and $\gamma_j:=\prod_{i=1}^j \beta_i$.
Hence we have $\gamma_s=D_0(\bfm,\bftau)$.
By Lemma~\ref{lem:gen-approx}, there exists $z_1\in F_{m_1-1,\tau_1}$ such that
$$
0<x-z_1\leq C_1 x^{\gamma_1}.
$$
Put $t_1=1+z_1\in F_{m_1,\tau_1}$.
Assume that for some $t_{j-1} \geq 1$, we have
$$
0<1+x-t_{j-1}\leq C_{j-1}x^{\gamma_{j-1}}.
$$
Write $\vep_j=(1+x-t_{j-1})/{t_{j-1}}$.
Lemma~\ref{lem:gen-approx} gives $z_j\in F_{m_j-1,\tau_j}$ such that
$$
0<\vep_j-z_j\leq C_j' {\vep_j}^{(\kappa_{\tau_j})^{m_j-1}}
 = \frac{C_j' (1+x-t_{j-1})^{ \beta_j }}{(t_{j-1})^{{ \beta_j }} }
 \leq   \frac{  C_j'(C_{j-1}x^{\gamma_{j-1}})^{ \beta_j }  }{(t_{j-1})^{ \beta_j } }     ,
$$
where $C_j'>0$ depends only on $m_j$ and $\tau_j$.
Write $C_j := C_j'  {(C_{j-1})}^{ \beta_j }$ and $t_j := t_{j-1}(1+z_j) \geq 1$.
Notice that $t_{j-1}, \beta_j \geq 1$, so that
$$
1+x-t_j
=t_{j-1}(\vep_j-z_j)
\leq C_j {t_{j-1}}^{1-{ \beta_j }}  x^{ \gamma_{j-1} \beta_j }
\leq C_j x^{ \gamma_j} .
$$
We get the desired result by iteration and $\gamma_s=D_0(\bfm,\bftau)$. 
The conclusion follows by taking $t=t_s$ and $C=C_s$.
\end{proof}

\begin{lem}\label{lem:gen-low-dim}
Let $s\geq 2$, $\bfm \in \Sigma^s$ and $\bftau \in (0,\infty)^s$. Then we have
$$
\lowdimB(\P_{\bfm, \bftau}) \geq 1-\frac{1}{D(\bfm,\bftau)}.
$$
\end{lem}

\begin{proof}
It suffices to prove that for each $0 \leq i \leq s$,
\begin{equation}\label{eq:gene-low-tmp}
\lowdimB( \P_{\bfm, \bftau} )\geq 1- 1/ D_i(\bfm,\bftau).
\end{equation}

By Lemma~\ref{lem:gen-gap-zero}, for every sufficiently small $\delta>0$, we have
$(1,1+\delta)\subset (\P_{\bfm, \bftau})_r. $
Here we take $r=C\delta^{D_0(\bfm,\bftau)} $.
\eqref{eq:min-dim-leb} then gives \eqref{eq:gene-low-tmp} for $i=0$.

By Lemma~\ref{lem:gen-gap-one}, for every sufficiently small $\delta>0$, we have
$(0,\delta)\subset (\P_{\bfm, \bftau})_r. $
Here we take $r=C\delta^{D_1(\bfm,\bftau)} $.
\eqref{eq:min-dim-leb} then gives \eqref{eq:gene-low-tmp} for $i=1$.

Fix $1 \leq j \leq s$ and write
$$\bfm'=(m_j, m_1 ,\ldots, m_{j-1}, m_{j+1} , \ldots,m_s ), \qquad
\bftau'=(\tau_j, \tau_1, \ldots, \tau_{j-1}, \tau_{j+1} , \ldots, \tau_s ).
$$
Since set multiplication is commutative, we have $\P_{\bfm, \bftau} = \P_{\bfm', \bftau'}$.
Applying the $i=1$ case to $\P_{\bfm', \bftau'}$, we have
$$
	\lowdimB( \P_{\bfm, \bftau} ) = \lowdimB( \P_{\bfm', \bftau'} ) \geq 1- 1/ D_1(\bfm',\bftau')
	= 1- 1/ D_j(\bfm,\bftau)  .
$$
Since $j$ was arbitrary, \eqref{eq:gene-low-tmp} holds for all $0 \leq i \leq s$.
Now the proof is complete.
\end{proof}

\subsection{Covering estimates and layer decompositions}

For any $m, a\in\Z^+$ and $\tau>0$, define
$$
U_{m, \tau ,a}:= 
\Big\{ a^{-\tau} +  x_1 + \cdots + x_{m-1}   :x_j\in F^{(\tau)} , x_j \leq a^{-\tau}, 1 \leq j \leq m - 1 \Big\}.
$$
Fix $s \in \Z^+$.
For any $\bfm=(m_1,\ldots,m_s), \bfa=(a_1,\ldots,a_s)  \in\Sigma^s$ and
$\bftau=(\tau_1,\ldots,\tau_s)\in(0,\infty)^s$, write
$$
L_{\bfm,\bftau,\bfa} :=  \prod_{i=1}^s U_{m_i, \tau_i ,a_i}.
$$

We next establish the layer decomposition for the generalized product sets.

\begin{lem}\label{lem:gen-layer}

Fix $s \in \Z^+$. 
For any $\bfm \in\Sigma^s$ and $\bftau \in (0,\infty)^s$, we have
$$
\P_{\bfm,\bftau} \setminus\{0\}
= \bigcup_{\bfa\in\Sigma^s} L_{\bfm,\bftau,\bfa}.
$$

\end{lem}

\begin{proof}
Fix $m\in\Z^+$ and $\tau>0$.
It is clear that $ 0 \notin U_{m,\tau,a}$ and $ U_{m,\tau,a} \subset F_{m, \tau} $ for every $a \in \Z^+$.
This implies that $\bigcup_{a\in\Z^+}   U_{m,\tau,a}\subset F_{m ,\tau } \setminus \{0\} $.
Conversely, let $x\in F_{m,\tau}\setminus\{0\}$.
Then $x=n_1^{-\tau}+\cdots+n_k^{-\tau}$ for some $1\leq k\leq m$ and $n_1\leq\cdots\leq n_k$.
Set $a=n_1$.
Then we have $x\in  U_{m,\tau,a}$. Since $x$ was arbitrary, $F_{m ,\tau } \setminus \{0\} \subset \bigcup_{a\in\Z^+}   U_{m,\tau,a}$.
Hence
\begin{equation}\label{eq:single-div}
F_{m,\tau}\setminus\{0\}=\bigcup_{a\in\Z^+}   U_{m,\tau,a}.
\end{equation}

Fix $\bfm=(m_1,\ldots,m_s) \in\Sigma^s$ and $\bftau=(\tau_1,\ldots,\tau_s)\in(0,\infty)^s$.
By \eqref{eq:single-div}, we obtain
$$
\P_{\bfm,\bftau} \setminus\{0\}
= \prod_{i=1}^s\big(F_{m_i,\tau_i}\setminus\{0\}\big)
= \bigcup_{\bfa\in\Sigma^s} \prod_{i=1}^s   U_{m_i,\tau_i,a_i}
= \bigcup_{\bfa\in\Sigma^s} L_{\bfm,\bftau,\bfa}.
$$
Then the proof is complete.
\end{proof}

In the remainder of this subsection, we are going to give a suitable estimate for the layer $ L_{\bfm,\bftau,\bfa}$.
For this purpose, we first normalize its factors.
For any $a\in\Z^+$ and $\tau>0$, define
$$
H_{\tau,a} :=  \big\{ (a/n)^\tau :n\in\Z^+,\ n\geq a \big\} \cup \{0\} .
$$

Fix $s\in\Z^+$.
For each $\bfm \in \Sigma^s$, write $K(\bfm):=\sum_{i=1}^s(m_i-1)$.
Here we may assume that $K(\bfm) \geq1$.
Define a function $\Phi:[0,1]^{K(\bfm)}\to\R$ by
\begin{equation}\label{eq:def-phi-k}
\Phi( x_{1,1},\ldots ,x_{1,m_1-1}, \ldots, x_{s,1},\ldots ,x_{s,m_s-1} )
:= \prod_{i=1}^s\bigg(1+\sum_{j=1}^{m_i-1}x_{i,j}\bigg).
\end{equation}
Then $\Phi$ is Lipschitz, and $\lip(\Phi)$ depends only on $\bfm$.
Moreover, for any $\bftau\in(0,\infty)^s$,
\begin{equation}\label{eq:layer-dec}
\prod_{i=1}^s a_i^{\tau_i} L_{\bfm,\bftau,\bfa}=  \prod_{i=1}^s a_i^{\tau_i} \prod_{i=1}^sU_{m_i ,\tau_i,a_i}
=   \Phi \big( H_{\tau_1,a_1}^{m_1-1}\times\cdots\times H_{\tau_s,a_s}^{m_s-1} \big).
\end{equation}
Here $H_{\tau,a}^m=\{(x_1,\ldots,x_m):x_i\in H_{\tau,a},\ 1\leq i\leq m\}$.
We now study Cartesian products of the sets $H_{\tau,a}$.
Fix $s\in\Z^+$.
For any $\bfb=(b_1,\ldots,b_s)\in\Sigma^s$, and $\bfsig=(\sigma_1,\ldots,\sigma_s)\in(0,\infty)^s$, we write
$$
E(\bfsig,\bfb) := \big\{ (x_1,\ldots,x_s): x_j\in H_{\sigma_j, b_j},\ x_1\geq\cdots\geq x_s \big\}.
$$
For each $ 1 \leq j \leq s $, write
$$
R_{\bfsig, j}:=\prod_{h=1}^j (\kappa_{\sigma_h})^{-1} = \prod_{h=1}^j \frac{\sigma_h}{1+\sigma_h}.
$$

\begin{lem}\label{lem:gen-harm-order-cover}
Let $s\in\Z^+$, $\bfb=(b_1,\ldots,b_s)\in\Sigma^s$, and $\bfsig=(\sigma_1,\ldots,\sigma_s)\in(0,\infty)^s$.
If $\Psi:[0,1]^s\to\R$ is Lipschitz, then for any $0<\rho<1$,
$$
\N_\rho\big(\Psi(E(\bfsig,\bfb))\big)
\leq C\rho^{R_{\bfsig,s} - 1}\prod_{j=1}^s b_j^{R_{\bfsig,j}},
$$
where $C$ depends only on $\bfsig$ and $\lip(\Psi)$.

\end{lem}

\begin{proof}
In this proof, write $\lambda:=\max\{1,\lip(\Psi)\}$.
Fix $\bfsig,\bfb$. Then we write $R_j = R_{\bfsig,j}$ for short.
We argue by induction on $s$.

For $s=1$, we have $R_1= (\kappa_{\sigma_1})^{-1}= \sigma_1/(1+\sigma_1)$.
Write $M :=  {b_1}^{R_1}  \rho^{R_1-1}$.
Then
$$
\N_\rho(H_{\sigma_1,b_1})
\leq \card \Big \{ x \in H_{\sigma_1,b_1} : x > \frac{{b_1}^{\sigma_1}}{M^{\sigma_1}} \Big\}
 +  \N_\rho \Big( \Big[0, \frac{{b_1}^{\sigma_1}}{M^{\sigma_1}}  \Big] \Big)
\leq M  + \frac{ {b_1}^{\sigma_1}}{\rho M^{\sigma_1}} +1.
$$
Since $ {b_1}^{\sigma_1} / \rho =M^{1+\sigma_1} $, we have $\N_\rho(H_{\sigma_1,b_1}) \leq 3M$.
Therefore
$$
\N_\rho\big(\Psi(E(\sigma_1,b_1)) \big)
= \N_\rho\big(\Psi(  H_{\sigma_1,b_1}  )\big)
\leq \N_{ \rho/\lambda}(H_{\sigma_1,b_1})
\leq  \frac{ 3 b_1^{R_1} }{ (\rho/\lambda)^{1-R_1} }
=  3\lambda^{1-R_1}\rho^{R_1-1} {b_1}^{R_1} .
$$
Taking $C_1:=3\lambda^{1-R_1}$ proves the assertion for $s=1$.

Fix $s\geq 2$ and assume the assertion holds for $s-1$.
Write
$\bfb'=(b_2,\ldots,b_s)$ and $\bfsig'=(\sigma_2,\ldots,\sigma_s)$.
For each $2 \leq j \leq s$, we set
$$
R_{j-1}' := R_{\bfsig',j-1} = \prod_{h=2}^j\frac{\sigma_h}{1+\sigma_h} = \frac{R_j}{R_1} .
$$
By the induction hypothesis, there exists a constant $C_{s-1}>0$ depending only on $\bfsig'$ and $\lambda $, such that for every Lipschitz map
$\Phi : [0,1]^{s-1} \to \R$ such that $\lip(\Phi)\leq \lambda $ and $ 0<\rho<1$,
\begin{equation}\label{eq:generalized-cover-ih}
\N_\rho\big( \Phi ( E(\bfsig',\bfb') ) \big)
\leq  C_{s-1} \rho^{ R_{s-1}'-1} \prod_{j=2}^s b_j^{ R_{j-1}' }
= C_{s-1} \rho^{(R_s / R_1 -1 )} \prod_{j=2}^s b_j^{ ( R_j / R_1 ) }.
\end{equation}

Write $B=\prod_{j=1}^s b_j^{R_j}$ and $T= \rho^{R_s}B$.
Since $B\geq1$ and $ R_s \in (0,1)$, we have $T / \rho  \geq 1$.
Write $C_s := C_{s-1}+s \lambda+1$.

First suppose that $T<1$.
For every $\bfx =(x_1, \ldots , x_s) \in E(\bfsig, \bfb ) $,
we have $\bfx \in [0, x_1]^s$ and $\bfx \in \{ x_1 \} \times E( \bfsig' , \bfb')$.
Set $ \Omega= \{ t \in H_{\sigma_1,b_1}: t >T\}$. 
Then
$$
E(\bfsig,\bfb) \subset [0,T]^s \cup \big( \Omega \times E(\bfsig' , \bfb') \big) .
$$
For each $t \in H_{\sigma_1,b_1}$, define $\Psi_t : [0,1]^{s-1} \to \R $ by
$$
\Psi_t (x_2 ,\ldots, x_s) := \Psi( t, x_2,\ldots,x_s).
$$
Then $\lip(\Psi_t )\leq\lambda$, and hence
\eqref{eq:generalized-cover-ih} gives
$$
\N_\rho\big( \Psi \big( \Omega \times E( \bfsig' , \bfb' ) \big) \big)
 \leq \sum_{t \in \Omega} \N_\rho\big(\Psi_t \big( E(\bfsig' , \bfb') \big) \big)
\leq \card (\Omega)  C_{s-1} \rho^{(R_s / R_1 -1)} \prod_{j=2}^s b_j^{(R_j / R_1)}.
$$
Since $ R_1 =1 - R_1/ \sigma_1$, we obtain
\begin{align*}
b_1 B^{-1 / \sigma_1}
\prod_{j=2}^s b_j^{ R_j / R_1}
= b_1^{1- R_1 / \sigma_1} \prod_{j=2}^s b_j^{ R_j / R_1 - R_j / \sigma_1}
= \prod_{j=1}^s b_j^{R_j}
= B.
\end{align*}
Combining with the fact that $\card ( \Omega ) \leq b_1T^{- 1 / \sigma_1 } =   \rho^{- R_s / \sigma_1 }  b_1 B^{- 1 / \sigma_1 }  $, we have
$$
\N_\rho\big( \Psi \big( \Omega \times E(\bfsig' , \bfb') \big) \big)
\leq  C_{s-1}   \rho^{(R_s / R_1 -1- R_s/ \sigma_1 )} B
=  C_{s-1}   \rho^{ R_s  -1 } B  .
$$
On the other hand, $\Psi([0,T]^s)$ can be covered by an interval of length $s T \lambda$.
Therefore
$$
\N_\rho\big(\Psi(E(\bfsig,\bfb))\big)
\leq C_{s-1} \rho^{R_s-1}B +s\lambda  \rho^{R_s}B / \rho+1
\leq C_s \rho^{R_s-1}B.
$$

If $T\geq1$, then $E(\bfsig,\bfb)\subset[0,1]^s$. Therefore
$$
\N_\rho\big(\Psi(E(\bfsig,\bfb))\big)
\leq s\lambda \rho^{-1}+1
\leq (s\lambda+1) T \rho^{-1}
\leq C_s \rho^{R_s -1 }B .
$$
This completes the induction.
\end{proof}

\begin{lem}\label{lem:gen-harm-cover-all-lip}
Let $s\in\Z^+$, $\bfb=(b_1,\ldots,b_s)\in\Sigma^s$, and $\bfsig=(\sigma_1,\ldots,\sigma_s)\in(0,\infty)^s$
such that $b_1^{\sigma_1}\geq\cdots\geq b_s^{\sigma_s}$.
Then for every Lipschitz map $\Phi:[0,1]^s\to\R$ and $0<\rho<1$, we have
$$
	\N_\rho \Big( \Phi \big(H_{\sigma_1,b_1}\times\cdots\times H_{\sigma_s,b_s} \big) \Big)
	\leq C \rho^{R_{\bfsig,s}-1} \prod_{j=1}^s b_j^{R_{\bfsig,j}},
$$
where the constant $C>0$ depends only on $\bfsig$ and $\lip(\Phi)$.
\end{lem}

\begin{proof}
Let $\Pi_s$ denote the set of permutations on $\{1,\ldots,s\}$.
For each permutation $\pi\in\Pi_s$, we write $\bfb_\pi=(b_{\pi(1)},\ldots,b_{\pi(s)})$ and $\bfsig_\pi=(\sigma_{\pi(1)},\ldots,\sigma_{\pi(s)})$.
Let $\Psi_\pi=\Phi\circ\pi^{-1}$. Then we have $\lip(\Psi_\pi)=\lip(\Phi)$.
Moreover,
$$
H_{\sigma_1,b_1}\times\cdots\times H_{\sigma_s,b_s} = \bigcup_{\pi\in\Pi_s}\pi^{-1} \big( E(\bfsig_\pi,\bfb_\pi) \big) .
$$

Define
$$
	g(\bfsig, \bfb, \pi) := \prod_{j=1}^s b_{\pi(j)}^{R_{\bfsig_\pi,j}}
$$
We prove by induction on $s$ that, for every $\pi\in\Pi_s$,
\begin{equation}\label{eq:per-sort}
g(\bfsig, \bfb, \pi) \leq g(\bfsig, \bfb, \id_s).
\end{equation}
Here $\id_s$ represents the identity permutation in $\Pi_s$.
The claim is trivial for $s=1$.
Now we suppose that it holds for $s-1$.

Consider $\pi =( \pi(1) , \ldots, \pi(s) ) \in \Pi_s$ with $\pi (i)=s$.
If $i=s$, the claim follows directly from the induction hypothesis. Thus we may assume that $1\leq i \leq s-1$.
Write
$$\pi' :=  \big( \pi(1) , \ldots  ,\pi(i-1), \pi(i+1), \ldots, \pi(s) , \pi(i) \big)  \in \Pi_s .$$
In this proof, for each $1 \leq j \leq s$, we write
$$
R_j:=\prod_{h=1}^j \frac{\sigma_{\pi'(h)}}{1+\sigma_{\pi'(h)}}.
$$
Set $R_0=1$ by convention. Then 
$$
R_{i-1}-R_{s-1}
=\sum_{j=i}^{s-1} (R_{j-1}-R_j )
=\sum_{j=i}^{s-1} R_j \bigg(\frac{1 +\sigma_{\pi'(j)} }{\sigma_{\pi'(j)} } -1  \bigg)
=\sum_{j=i}^{s-1}\frac{R_j}{\sigma_{\pi'(j)}}.
$$
Since $b_s= b_{\pi'(s)}=b_{\pi(i)}$, by definition,
$$
\frac{g(\bfsig, \bfb, \pi') }{ g(\bfsig, \bfb, \pi) }
= (b_s)^{ R_{s} - \frac{\sigma_s R_{i-1}}{1+\sigma_s}}
\prod_{j=i}^{s-1} \big(b_{\pi'(j)}\big)^{ R_{j}  -\frac{\sigma_s R_{j}  }{1+\sigma_s} }
=  \big({b_s}^{\sigma_s}\big)^{   \frac{R_{s-1}-R_{i-1} }{1+\sigma_s}}
 \prod_{j=i}^{s-1} \big(  b_{\pi'(j)} \big)^{ \frac{ R_j  }{ 1+\sigma_s} } .
$$
Recall that $(b_{\pi'(j)})^{\sigma_{\pi'(j)} }\geq b_s^{\sigma_s}$ for each $1 \leq j \leq s$.
Let $q=b_s^{\sigma_s}$. 
Therefore
$$
 \bigg( \frac{g(\bfsig, \bfb, \pi') }{ g(\bfsig, \bfb, \pi) } \bigg)^{1 +\sigma_s}
=  q^{  R_{s-1}-R_{i-1} }
 \prod_{j=i}^{s-1} \big(  (b_{\pi'(j)})^{\sigma_{\pi'(j)} }  \big)^{ \frac{ R_j  }{\sigma_{\pi'(j)} } }
 \geq   q^{  R_{s-1}-R_{i-1} }
 \prod_{j=i}^{s-1} q^{ \frac{ R_j  }{\sigma_{\pi'(j)} } } =1.
$$
This implies that $g(\bfsig, \bfb, \pi') \geq g(\bfsig, \bfb, \pi) $.
Since $\pi'$ and $\id_s$ share the same last entry, by the induction hypothesis, we have
$$
	g(\bfsig, \bfb, \id_s) \geq g(\bfsig, \bfb, \pi') \geq g(\bfsig, \bfb, \pi) .
$$
Let $C_\pi$ be the constant in Lemma~\ref{lem:gen-harm-order-cover} depending on $\bfsig_\pi$ and $\lip(\Phi)=\lip(\Psi_\pi)$.
Here we write $C':=\max_{\pi\in\Pi_s}C_\pi$. Then 
$$
\N_\rho\big(\Phi(\pi^{-1}E(\bfsig_\pi,\bfb_\pi))\big)
= \N_\rho\big(\Psi_\pi(E(\bfsig_\pi,\bfb_\pi))\big)
\leq C'\rho^{R_s-1} g(\bfsig, \bfb, \pi)
\leq C' \rho^{R_s-1} g(\bfsig, \bfb, \id_s) .
$$
The last inequality follows from \eqref{eq:per-sort}.
Summing over $\pi\in\Pi_s$, we have
$$
	\N_\rho \Big( \Phi \big(H_{\sigma_1,b_1}\times\cdots\times H_{\sigma_s,b_s} \big) \Big)
	\leq \sum_{\pi \in \Pi_s} \N_\rho\big(\Phi(\pi^{-1}E(\bfsig_\pi,\bfb_\pi))\big)
	\leq s! C' \rho^{R_{\bfsig,s}-1} \prod_{j=1}^s b_j^{R_{\bfsig,j}}.
$$
Writing $C :=s! C'$ gives the desired estimate.
\end{proof}

Applying Lemma~\ref{lem:gen-harm-cover-all-lip} to \eqref{eq:layer-dec} gives the following covering estimate.

\begin{lem}\label{lem:gen-block-cover}
Let $s\in\Z^+$, $\bfm=(m_1,\ldots,m_s),\bfa=(a_1,\ldots,a_s)\in\Sigma^s$, and $\bftau=(\tau_1,\ldots,\tau_s)\in(0,\infty)^s$.
Then there exists a constant $C>0$, depending only on $\bfm$ and $\bftau$, 
such that for every $0<\rho < \prod_{i=1}^s {a_i}^{ -\tau_i}$,
\begin{equation}\label{eq:generalized-layer-cover-new}
\N_\rho(L_{\bfm,\bftau,\bfa})
\leq C \rho^{ D_0(\bfm,\bftau)^{-1}-1 }
\prod_{i=1}^s (a_i^{\tau_i})^{D_0(\bfm,\bftau)^{-1} - {(\kappa_{\tau_i})}^{1-m_i}} ,
\end{equation}

\end{lem}

\begin{proof}
Write $K = K(\bfm) =\sum_{i=1}^s(m_i-1)$.
If $K=0$, then we have $m_i=1$ for every $1\leq i\leq s$.
Thus $D_0(\bfm,\bftau)=1$ and $L_{\bfm,\bftau,\bfa}$ is a singleton.
By taking $C=1$, the desired estimate is immediate.
Hence we may assume that $K\geq1$.

Let $\Phi$ be the function defined in \eqref{eq:def-phi-k}.
We remark that $\Phi$ is Lipschitz, and $\lip(\Phi)$ depends only on $\bfm$.
Let $A:=\prod_{i=1}^s a_i^{\tau_i}$.
Then we have $\rho A<1$.

Without loss of generality, we assume that $a_1^{\tau_1}\geq\cdots\geq a_s^{\tau_s}$.
By Lemma~\ref{lem:gen-harm-cover-all-lip} and \eqref{eq:layer-dec},
$$
\N_\rho(L_{\bfm,\bftau,\bfa})
= \N_{\rho A }\Big(  \Phi \big( H_{\tau_1,a_1}^{m_1-1}\times\cdots\times H_{\tau_s,a_s}^{m_s-1} \big) \Big)
\leq C(\rho A)^{D_0(\bfm,\bftau)^{-1} -1 } \prod_{i=1}^s a_i^{ \tau_i(1-\kappa_{\tau_i}^{1-m_i}) }  .
$$
Here the exponent of $a_i$ is obtained by
$$
 \sum_{j=1}^{m_i-1}\kappa_{\tau_i}^{-j} \prod\limits_{h=1}^{i-1}\kappa_{\tau_h}^{1-m_h}
\leq \sum_{j=1}^{m_i-1} \kappa_{\tau_i}^{-j}
= \frac{ \kappa_{\tau_i}^{-1} (1- \kappa_{\tau_i}^{-(m_i-1)} ) }{1 - \kappa_{\tau_i}^{-1}}
=\tau_i(1-\kappa_{\tau_i}^{1-m_i}).
$$
Therefore,
$$
\begin{aligned}
\N_\rho(L_{\bfm,\bftau,\bfa})
\leq
C(\rho A)^{D_0(\bfm,\bftau)^{-1} -1}
\prod_{i=1}^s a_i^{\tau_i(1-\kappa_{\tau_i}^{1-m_i})}
=
C\rho^{D_0(\bfm,\bftau)^{-1} -1}
\prod_{i=1}^s
\big(a_i^{\tau_i}\big)^{D_0(\bfm,\bftau)^{-1} - (\kappa_{\tau_i})^{1-m_i} },
\end{aligned}
$$
which proves \eqref{eq:generalized-layer-cover-new} for the case $a_1^{\tau_1}\geq\cdots\geq a_s^{\tau_s}$.
Since the estimate in \eqref{eq:generalized-layer-cover-new} is invariant under simultaneous reordering of the parameters, it holds for arbitrary $\bfa$ and $\bftau$.
\end{proof}

\subsection{The upper bound}

We are going to combine the preceding covering estimates with the layer decomposition.
Fix $s,n \in \Z^+$. For any $\bftau \in (0, \infty)^s$, write 
$$\Lambda_{\bftau,n} := \{ \bfa\in\Sigma^s :  \prod_{i=1}^s  a_i^{\tau_i}   \leq n\}.$$
To sum these layer estimates, we introduce the following index bound.

\begin{lem}\label{lem:gen-index-sum}
Let $s\in\Z^+$, $\bftau=(\tau_1,\ldots,\tau_s)\in(0,\infty)^s$, and
$\bfeta=(\eta_1,\ldots,\eta_s)\in (-\infty,0 ]^s$.
Assume that $\beta \geq 0$ and $1/\tau_i + \eta_i\leq\beta$ for every $1\leq i\leq s$.
Then for every $ n \geq2$,
$$
\sum_{ \bfa \in \Lambda_{\bftau,n} }
\prod_{i=1}^s a_i^{\tau_i \eta_i}
\leq C n^\beta(\log n)^{s},
$$
where $C>0$ depends only on $\bftau$.
\end{lem}

\begin{proof}
Fix $\bftau$.
Write $C_{\bftau} := \prod_{i=1}^s (2+ 1/ \tau_i)$.
Fix $\tau>0$ and $\eta\leq 0$ such that $1/\tau + \eta\leq\beta$. 
Then for every $X\geq 2$, we have
\begin{equation}\label{eq:generalized-one-index-sum-0}
\sum_{a^{\tau}  \leq X} a^{\tau(\eta - \beta)}
\leq \sum_{a^{\tau}  \leq X} a^{-1}
\leq 1 + \frac{\log X}{ \tau}
\leq \frac{2\tau +1 }{\tau} \log X.
\end{equation}

Since $\beta\geq0$, we have
$$
\sum_{ \bfa \in \Lambda_{\bftau,n} } \prod_{i=1}^s a_i^{\tau_i\eta_i}
=\sum_{ \bfa \in \Lambda_{\bftau,n} } \bigg(\prod_{i=1}^s a_i^{\tau_i}\bigg)^\beta
\prod_{i=1}^s a_i^{\tau_i(\eta_i-\beta)} 
\leq n^\beta \prod_{i=1}^s \sum_{a_i^{\tau_i}\leq n } a_i^{\tau_i(\eta_i-\beta)}
$$
By the assumption that $1/\tau_i + \eta_i\leq\beta$ and \eqref{eq:generalized-one-index-sum-0},
we have
$$
\prod_{i=1}^s \sum_{a_i^{\tau_i}\leq n } a_i^{\tau_i(\eta_i-\beta)}
\leq \prod_{i=1}^s \Big( \frac{2\tau_i +1 }{\tau_i} \log n \Big)
= C_\bftau (\log n)^s.
$$
Taking $C=C_{\bftau}$ completes the proof.
\end{proof}

\begin{lem}\label{lem:gen-layer-cover}
Let $s\in\Z^+$, $\bfm=(m_1,\ldots,m_s) \in\Sigma^s$, and $\bftau=(\tau_1,\ldots,\tau_s)\in(0,\infty)^s$.
Then there exists a constant $C>0$ depending only on $\bfm$ and $\bftau$ such that for every $n\geq2$,
$$
\sum_{\bfa\in \Lambda_{\bftau,n}}
\N_{n^{-D(\bfm,\bftau)}}(L_{\bfm,\bftau,\bfa})
\leq C n^{D(\bfm,\bftau)-1}(\log n)^s,
$$

\end{lem}

\begin{proof}
Write $P:=D_0(\bfm,\bftau)^{-1}$, $D:=D(\bfm,\bftau)$ and $\rho:=n^{-D}$.
For each $1\leq i\leq s$, we have 
$$
	D_0(\bfm,\bftau)^{-1} =\prod_{j=1}^s (\kappa_{\tau_j})^{1 - m_j} \leq (\kappa_{\tau_i})^{1-m_i}.
$$
Set $ \eta_i:= D_0(\bfm,\bftau)^{-1} -\kappa_{\tau_i}^{1-m_i}  \leq 0 $.
Hence $\bfeta = ( \eta_1,\ldots,\eta_s  ) \in(-\infty,0]^s$.
By \eqref{eq:generalized-D-short},
$$
D_i(\bfm,\bftau)
=1+\big(\kappa_{\tau_i}-\kappa_{\tau_i}^{1-m_i}\big)
D_0(\bfm,\bftau).
$$
Using $\kappa_{\tau_i}-1=1/\tau_i$ and $D$ defined in \eqref{eq:generalized-D}, we obtain
\begin{equation}\label{eq:def-tmp-eta-pri}
\frac{1}{\tau_i}+\eta_i
=P+\frac{1}{\tau_i}-\kappa_{\tau_i}^{1-m_i}
=PD_i(\bfm,\bftau)-1
\leq PD-1.
\end{equation}

Write $\beta:=PD-1$.
Then $\beta\geq 0$ follows from $D\geq D_0(\bfm,\bftau)=P^{-1}$,
 while
$1 / \tau_i +\eta_i\leq\beta$ follows from \eqref{eq:def-tmp-eta-pri}.

For every $\bfa\in\Lambda_{\bftau,n}$, we have
$$
\rho\prod_{i=1}^s a_i^{\tau_i} \leq n^{-D} n = n^{1-D} <1.
$$
Then Lemma~\ref{lem:gen-block-cover} gives
$$
\N_\rho(L_{\bfm,\bftau,\bfa})
\leq C \rho^{P-1}
\prod_{i=1}^s a_i^{\tau_i\eta_i}.
$$
By Lemma~\ref{lem:gen-index-sum},
$$
\sum_{\bfa\in\Lambda_{\bftau,n}} \N_\rho(L_{\bfm,\bftau,\bfa})
\leq C\rho^{P-1} \sum_{\bfa \in\Lambda_{\bftau,n} } \prod_{i=1}^s a_i^{\tau_i\eta_i}
\leq
C\rho^{P-1}n^{PD-1}(\log n)^s.
$$
By definition, we have $\rho^{P-1}n^{PD-1}=n^{D-1}$,
which completes the proof.
\end{proof}

\begin{lem}\label{lem:gen-up-dim}
Let $\bfm=(m_1,\ldots,m_s)\in\Sigma^s$ and $\bftau=(\tau_1,\ldots,\tau_s)\in(0,\infty)^s$. Then
$$
\updimB(\P_{\bfm,\bftau}) \leq 1-\frac{1}{D(\bfm,\bftau)}.
$$
\end{lem}

\begin{proof}

By Lemma~\ref{lem:gen-layer}, we have
$$
\begin{aligned}
\P_{\bfm,\bftau}
=
\{0\}\cup
\bigcup_{\bfa\in\Sigma^s}
L_{\bfm,\bftau,\bfa}.
\end{aligned}
$$

Fix $\bfm,\bftau$ and an integer $n\geq 3$ and set $\rho:=n^{-D(\bfm,\bftau)}$.
For each $\bfa=(a_1,\ldots,a_s)\in\Sigma^s \setminus \Lambda_{\bftau,n} $, 
we have $ {a_i}^{\tau_i} U_{m_i,\tau_i,a_i}\subset[1,m_i]$ for each $1 \leq i \leq s$.
Then we have
$$
\bigcup_{\bfa\in \Sigma^s \setminus \Lambda_{\bftau,n}} L_{\bfm,\bftau,\bfa} 
\subset \bigcup_{\bfa\in \Sigma^s \setminus \Lambda_{\bftau,n}}
 \Big[  \prod_{i=1}^s  \frac{1}{ a_i^{\tau_i} } , \prod_{i=1}^s  \frac{m_i}{ a_i^{\tau_i} } \Big]
\subset \Big[ 0,\frac{1}{n}  \prod_{i=1}^s m_i  \Big].
$$
Write $ M :=\prod_{i=1}^s m_i  $. 
Then
\begin{equation}\label{eq:generalized-tail-cover-new}
\N_\rho\Big( \{0\} \cup \bigcup_{ \bfa\in\Sigma^s \setminus\Lambda_{\bftau,n}  } L_{\bfm,\bftau,\bfa} \Big)
\leq \N_\rho\Big( \Big[ 0,\frac{M}{n} \Big] \Big)
\leq 2M n^{D(\bfm,\bftau)-1} (\log n)^s.
\end{equation}

By Lemma~\ref{lem:gen-layer-cover}, we obtain
\begin{equation}\label{eq:generalized-main-cover-new}
\sum_{ \bfa\in \Lambda_{\bftau,n}} \N_\rho(L_{\bfm,\bftau,\bfa})
\leq C n^{D(\bfm,\bftau)-1}(\log n)^s.
\end{equation}  
Combining \eqref{eq:generalized-tail-cover-new} and
\eqref{eq:generalized-main-cover-new}, we have
$$
\N_{n^{-D(\bfm,\bftau)}}(\P_{\bfm,\bftau})
\leq (2M +C) n^{D(\bfm,\bftau)-1}(\log n)^s.
$$
For every sufficiently small $r>0$, choose $n\geq2$ such that
$
n^{-D(\bfm,\bftau)} \leq r < (n-1)^{-D(\bfm,\bftau)}.
$
Then
$$
\N_r(\P_{\bfm,\bftau})
\leq
\N_{n^{-D(\bfm,\bftau)}}(\P_{\bfm,\bftau}),
$$
and hence by \eqref{eq:min-dim-cover}
$$
\updimB(\P_{\bfm,\bftau})
\leq
\varlimsup_{n\to\infty}
\frac{ \log\big( (2M+C) n^{D(\bfm,\bftau)-1}(\log n)^s\big)}{D(\bfm,\bftau)\log(n-1)}
=  \frac{D(\bfm,\bftau) - 1}{D(\bfm,\bftau)}.
$$
This completes the proof.
\end{proof}

\begin{proof}[Proof of Theorem~\ref{thm:gen}]
If $s=1$, the conclusion follows from 
$$\dimB(F_{m_1,\tau_1})=1-\kappa_{\tau_1}^{-m_1}. $$
Hence we assume $s\geq2$.
The lower bound is given by Lemma~\ref{lem:gen-low-dim},
while the upper bound follows from Lemma~\ref{lem:gen-up-dim}.
\end{proof}

We conclude this section with several useful special cases of Theorem~\ref{thm:gen}.

\begin{rem}\label{rem:gen-fast-1}
For every $p,q\geq2$ and $\tau,\sigma>0$, we have
$$
\dimB\big(F_{p,\tau }F_{q,\sigma }\big) =1-
\frac{1}{\max\big\{
	\kappa_\tau^{p-1} \kappa_\sigma^{q-1},  1+(\kappa_\tau^p-1)\kappa_\sigma^{q-1} ,  1+(\kappa_\sigma^q-1)\kappa_\tau^{p-1}
	\big\}}.
$$

\end{rem}

\begin{rem}\label{rem:gen-fast-2}
Fix $s \in \Z^+$ and $\tau>0$.
Let $\bfm=(m_1,\ldots,m_s)\in\Sigma^s$ and $\bftau:=(\tau,\ldots,\tau)$,
$$
K:=\sum_{i=1}^s(m_i-1),
\qquad
m_*:=\max_{1\leq i\leq s}m_i.
$$
Then for each $1\leq i\leq s$,
\begin{equation}\label{eq-rem-tau-di}
D_i(\bfm,\bftau)
= 1+(\kappa_\tau^{m_i}-1)\kappa_\tau^{K-m_i+1}
= 1+\kappa_\tau^{K+1}-\kappa_\tau^{K-m_i+1}.
\end{equation}
Indeed, we have $ D_0(\bfm,\bftau)=\kappa_\tau^K$. 
Since $\kappa_\tau>1$, $\max_{1 \leq i \leq s}D_i$ is achieved when $m_i=m_*$. 
Hence Theorem~\ref{thm:gen} gives
$$
\dimB(\P_{\bfm,\bftau}) = 1-\frac{1}{\max \{ \kappa_\tau^K,\, 1+(\kappa_\tau^{m_*}-1) \kappa_\tau^{K-m_*+1} \}}.
$$

In particular, if $\bfsig = (1, \ldots , 1) \in (0,\infty)^s$, then we have $ \kappa_{\sigma_i}=2 $ for each $1 \leq i \leq s$. 
Then we have $D_0(\bfm, \bfsig) = 2^K$.
By \eqref{eq-rem-tau-di},
$$
D_i(\bfm,\bfsig)
= 1+2^{K+1}-2^{K-m_i+1} \geq 1+2^K \geq D_0(\bfm, \bfsig).
$$
Hence we obtain
$$
	D(\bfm,\bfsig) =\max_{0 \leq i\leq s} D_i(\bfm,\bfsig) 
	=  1 + (2 - 2^{1-m_*}) 2^K
	=  1+(2^{m_*}-1)2^{K- m_* +1}.
$$
This proves Theorem~\ref{thm:prod}.

\end{rem}

\section{Taylor sets}

In this section, we prove the lower bound for $\T_{k,m}$ and then establish the upper bound for $\T_{2,m}$.

\subsection{The lower bound}

For every $m,s \in \Z^+$, write $(F_m)^s=\prod_{j=1}^s F_m$.
Then 
$$
	\T_{k,m} = \sum_{i=1}^k (F_m)^i.
$$

For products of $F_m$, the constant in Lemma~\ref{lem:gen-gap-zero} can be taken to be $1$.

\begin{lem}\label{lem:true-taylor-product-gap}
Let $s,m\in\Z^+$ with $m\geq2$. For every $x\in(0,1)$, there exists $u\in(F_m)^s$ satisfying
$$
0\leq x +1 -u\leq x^{2^{s(m-1)}}.
$$
\end{lem}

\begin{proof}
Without loss of generality, we assume that $1 + x \notin (F_m)^s$.
For any $z\in(0,1)$, choose $n\in\Z^+$ such that $ 1/ n<z\leq 1/(n-1)$. Then
\begin{equation}\label{eq:tay-low-fast-tmp-1}
0<z-\frac{1}{n}\leq\frac{z}{n}<z^2.
\end{equation}
By iterating \eqref{eq:tay-low-fast-tmp-1}, there exists $v\in F_{m-1}$ such that
$
0<z-v\leq z^{2^{m-1}}.
$

Set $u_0:=1$. For $1\leq i \leq s$, suppose that $1\leq u_{i-1}<1+x$. 
Let $z_i = (x+1)/u_{i-1}-1 >0$.
Then there exists $v_i \in F_{m-1}$ such that $ 0<z_i -v_i \leq z_i^{2^{m-1}}$.
Set $u_i := u_{i-1}(1+v_i) $. 
Then $u_i \in (F_m)^i $ and $1 \leq u_i <1+x $.
 Since $u_{i-1}\geq1$, we obtain
$$
0<x +1 -u_i
\leq u_{i-1}^{1-2^{m-1}} (x+1 -u_{i-1})^{2^{m-1}}
\leq (x+1-u_{i-1})^{2^{m-1}}.
$$
By induction, we obtain
$$
 \log( x +1 -u_s ) \leq 2^{m-1} \log (x+1-u_{s-1}) \leq  \cdots \leq 2^{s(m-1)} \log(x+1-u_{0}) = 2^{s(m-1)}\log x.
$$
Taking $u:=u_s$ completes the proof.
\end{proof}

\begin{lem}\label{lem:true-taylor-approx}
Let $k,m\in\Z^+$ with $m\geq2$.
Then for every $ x \in (0, 1/4)$, there exists
$y \in \T_{k,m}$ satisfying
$$
0<x-y \leq (2x)^{2^{km}}.
$$
\end{lem}

\begin{proof}
Fix $m\geq2$. We prove the lemma by induction on $k$.
For $k=1$, iterating \eqref{eq:tay-low-fast-tmp-1} $m$ times gives $a\in F_m=\T_{1,m}$ such that
$$
0<x-a\leq x^{2^m}\leq(2x)^{2^m}.
$$

Now suppose that $k\geq2$ and that the assertion holds for $k-1$.
Fix $x\in(0, 1 /4 )$.
Choose $n\in\Z^+$ such that $ 1 / n < x \leq 1 /(n-1) $. 
Set $z := x -1/(n+1)$.
Since $n\geq5$, we have
$$
\frac{x^2}{2} < \frac{x}{n+1} <z \leq  \frac{2x}{n+1} \leq 2x^2<\frac{1}{4}.
$$
Applying the construction for $k-1$ to $z$, there exists $b\in\T_{k-1,m}$ such that
$$
0<z-b \leq (2 z)^{2^{(k-1)m}}.
$$
Since $0<z-b<1$, iterating \eqref{eq:tay-low-fast-tmp-1} $m-1$ times gives $v\in F_{m-1}$ such that
$$
0<z-b-v\leq(z-b)^{2^{m-1}}.
$$

Notice that $v+ 1/(n+1)\in F_m$ and $1\in F_m$.
Set $a:= v+ 1/(n+1)+b$.
Then we have
$
a\in(F_m)^k+\T_{k-1,m}=\T_{k,m}.
$
Moreover,
$$
0<x-a
=z-b-v
\leq(z-b)^{2^{m-1}}
\leq(2z)^{2^{km-1}}
\leq(2x)^{2^{km}}.
$$
This completes the induction.
\end{proof}

Based on Lemma~\ref{lem:true-taylor-approx}, we obtain the following result.

\begin{lem}\label{lem:true-taylor-gap}
Let $k,m\in\Z^+$ with $m\geq2$. There exists a constant $C>0$, depending only on $k$ and $m$, such that for every $x\in(0,1/4)$, there exists $y\in\T_{k,m}$ satisfying
$$
0\leq x-y\leq Cx^{d_{k,m}}.
$$
\end{lem}

\begin{proof}
Fix $m\geq2$. 
For each $k\in\Z^+$, write
$$
p_k:=2^{ k(k-1) (m-1)/ 2}.
$$
Then we have $p_1=1$ and $p_k=2^{(k-1)(m-1)}p_{k-1}$ for every $ k \geq 2$.

We claim that for every $x\in(0,1/4)$, there exists $a\in\T_{k,m}$ with $0<x-a\leq(2x)^{2^{km}}$ such that, for every $0\leq\vep\leq(x/2)^k$, there exists $y\in\T_{k,m}$ satisfying
\begin{equation}\label{eq:taylor-gap-induction}
0\leq a+\vep-y
\leq \Big(\frac{x}{2}\Big)^{d_{k,m}-2^{km}p_k}\vep^{p_k}.
\end{equation}
We prove by induction on $k$. For $k=1$, $d_{1,m}=2^m$.
Choose $a$ by Lemma~\ref{lem:true-taylor-approx} and take $y=a$. Then the assertion follows.

Let $k\geq2$. 
Choose $n\in\Z^+$ such that $1/n<x\leq 1/(n-1)$.
Write $z:=x-1/(n+1)$. Hence $x^2 / 2 <z<2x^2<1 /4 $.
By the claim at $z$, there exists $a'\in\T_{k-1,m}$.
Iterating \eqref{eq:tay-low-fast-tmp-1} $m-1$ times, there exists $b\in 1/(n+1)+F_{m-1}$ such that
$$
0<x-a'-b\leq(z-a')^{2^{m-1}} \leq(2z)^{2^{km-1}} \leq(2x)^{2^{km}}.
$$
Then $a:=b+a' \in \T_{k,m}$ and $x/2<b<x$.
Fix $0\leq\vep\leq(x/2)^k$.
If $\vep=0$, take $y=a$. Otherwise, $0<\vep/b\leq(x/2)^{k-1}<1$.
By Lemma~\ref{lem:true-taylor-product-gap}, 
there exists $u / b \in (F_m)^{k-1}$ satisfying
$$
0\leq b+\vep-u
\leq \Big(\frac{x}{2} \Big)^{1-2^{(k-1)(m-1)}} \vep^{2^{(k-1)(m-1)}}.
$$
Then $u \in (F_m)^k$.
For simplicity, write $q := 2^{(k-1)(m-1)} \geq 2$.	
Since $\vep\leq (x/2)^k$, 
\begin{equation}\label{eq:taylor-gap-tmp-1}
0 \leq b+\vep-u
\leq \Big(\frac{x}{2}\Big)^{1-q} \vep^{q}
\leq \Big(\frac{x}{2}\Big)^{1-q} \Big(\frac{x}{2}\Big)^{kq}
\leq \Big(\frac{x}{2}\Big)^{2(k-1)}
\leq \Big(\frac{z}{2}\Big)^{k-1}.
\end{equation}
Notice that $a+\vep-u=a'+(b+\vep-u)$.
Applying the claim at $z$ with $b+\vep-u$, we obtain
$y\in u+\T_{k-1,m}\subset\T_{k,m}$ satisfying
$$
0\leq a+\vep-y
\leq \Big(\frac{z}{2}\Big)^{d_{k-1,m}-2^{(k-1)m}p_{k-1}} (b+\vep-u)^{p_{k-1}}.
$$

By \eqref{eq:true-taylor-D} and the definition of $p_k$, we obtain
$$
d_{k,m}-2^{km}p_k
= 2(d_{k-1,m}-2^{(k-1)m}p_{k-1}) +p_{k-1}-p_k.
$$
By~\eqref{eq:taylor-gap-tmp-1}, we obtain that
$$
(b+\vep-u)^{p_{k-1}}
\leq \Big( \Big(\frac{x}{2}\Big)^{1-q}\vep^q \Big)^{p_{k-1}}
= \Big(\frac{x}{2}\Big)^{p_{k-1}-p_k} \vep^{p_k}.
$$
Since $d_{k-1,m}-2^{(k-1)m}p_{k-1}\leq0$ and $z/2>(x/2)^2$,
$$
a+\vep-y
\leq
\Big(\frac{x}{2}\Big)^{2(d_{k-1,m}-2^{(k-1)m}p_{k-1})}
\Big(\frac{x}{2}\Big)^{p_{k-1}-p_k} \vep^{p_k}
=\Big(\frac{x}{2}\Big)^{d_{k,m}-2^{km}p_k} \vep^{p_k}.
$$
Finally, take $\vep=x-a$. Since $2^{km}\geq3k$ and $x<1/4$, we have
$$
0<\vep\leq(2x)^{2^{km}} \leq(2x)^{3k} \leq(x/2)^k.
$$
By \eqref{eq:taylor-gap-induction}, there exists $y\in\T_{k,m}$ such that
$$
0\leq x-y
\leq \Big(\frac{x}{2}\Big)^{d_{k,m}-2^{km}p_k} \big( (2x)^{2^{km}} \big)^{p_k}
= 2^{2^{km+1}p_k-d_{k,m}}x^{d_{k,m}}.
$$
This proves the lemma.
\end{proof}

\begin{lem}\label{lem:true-taylor-low-dim}
For every $k,m\in\Z^+$ with $m \geq 2$, we have
$$
\lowdimB(\T_{k,m})\geq1-\frac{1}{d_{k,m}}.
$$
\end{lem}

\begin{proof}
Let $C$ be the constant in Lemma~\ref{lem:true-taylor-gap}. For every sufficiently small $\delta>0$, set $r:=C\delta^{d_{k,m}}$. Then
$$
(0,\delta)\subset(\T_{k,m})_r,
\qquad
\leb((\T_{k,m})_r)\geq\delta.
$$
The desired estimate follows from \eqref{eq:min-dim-leb}.
\end{proof}

\subsection{The upper bound}

For $m ,a \in\Z^+$, define
$$
F_m(a):=
\Big\{
\sum_{i=1}^m x_i:
x_i\in F,\ x_i\leq \frac{1}{a}
\Big\}.
$$
We also set $F_0(a):=\{0\}$ by convention.

\begin{lem}\label{lem:s3-single-layer-cover}
Let $m\in\Z^+$ with $m\geq2$.
There exists a constant $C_m >0$ such that for every $0 \leq k \leq m$, $a,b\in\Z^+$ with $a\geq b$ and every $0<\rho<(ab)^{-1}$,
$$
\N_\rho\big( F_k (ab) + U_{m,1,a}U_{m,1,b} \big)
\leq C_m\rho^{2^{2-k-2m}-1} (a^{2^{2-k-2m}}b^{2^{1-k-m}})^{1-2^{m-1}}.
$$
\end{lem}

\begin{proof}
Notice that for every $a \in \Z^+$, we have
$$
	 a F_m ( a ) = \sum_{i=1}^m H_{1,a} =aU_{m+1,1,a}-1 .
$$
Define a function $\Phi$ on $[0,1]^{k+2m-2}$ by
$$
\Phi(u_1,\ldots,u_k,x_1,\ldots,x_{m-1},y_1,\ldots,y_{m-1})
:= \sum_{i=1}^k u_i + \Big(1+\sum_{i=1}^{m-1}x_i\Big) \Big(1+\sum_{i=1}^{m-1}y_i\Big).
$$
Then $\Phi$ is Lipschitz, and $\lip(\Phi) \leq\sqrt{k+2m^2(m-1)} \leq m^2$.
Moreover,
$$
 F_k (ab) + U_{m,1,a}U_{m,1,b} =
\frac{1}{ab} \Phi\Big( H_{1,ab}^{k}\times H_{1,a}^{m-1}\times H_{1,b}^{m-1} \Big).
$$
Since $ab\geq a\geq b$ and $ab\rho<1$, Lemma~\ref{lem:gen-harm-cover-all-lip} gives
$$
\begin{aligned}
&\N_\rho\big(F_k(ab)+U_{m,1,a}U_{m,1,b}\big)\\
\leq&
C_m(ab\rho)^{2^{-(k+2m-2)}-1}
(ab)^{1-2^{-k}}
a^{2^{-k}(1-2^{1-m})}
b^{2^{1-k-m}(1-2^{1-m})}\\
=&
C_m\rho^{2^{2-k-2m}-1}
(a^{2^{2-k-2m}}b^{2^{1-k-m}})^{1-2^{m-1}}.
\end{aligned}
$$
\end{proof}

The layer decomposition now yields the following covering estimate.

\begin{lem}\label{lem:s3-tail-cover}
Let $m \in\Z^+$ with $m\geq2$. Let $0 \leq k \leq m$.
There exists a constant $C>0$, depending on $m,k$ such that for every $h\in \Z^+$ and $0<\rho<h^{-1}$,
$$
\N_\rho\big(F_k(h) + F_mF_m\big)
\leq C \big(  \rho^{-1} h^{-1} + \rho^{2^{2-k-2m}-1} h^{1-(2^{m-1}-1)2^{2-k-2m}} \big).
$$
\end{lem}

\begin{proof}

For any $0<\tau<\sigma<1$ and $X \geq 1$, integral comparison gives
$$
\sum_{t \leq X: t \in \Z^+} t^{-\tau} \leq1+\int_1^X t^{-\tau}\,dt \leq\frac{X^{1-\tau}}{1-\tau} .
$$
Then we have
\begin{equation}\label{eq:s3-asymmetric-sum}
\sum_{ab\leq h : a\geq b} a^{-\tau}b^{-\sigma}
\leq
\sum_{b\leq h} b^{-\sigma}
\sum_{a\leq h/b}a^{-\tau} 
\leq \frac{h^{1-\tau}}{1-\tau} \sum_{b=1}^{\infty}b^{- \sigma -1 + \tau}
\leq \frac{h^{1-\tau}}{1-\tau} \Big(1+\frac{1}{\sigma-\tau}\Big).
\end{equation}

By Lemma~\ref{lem:gen-layer},
$$
F_mF_m
\subset
\Big[ 0,\frac{m^2}{h} \Big]
\cup
\bigcup_{ab\leq h} U_{m,1,a}U_{m,1,b}.
$$
Hence
$$
F_k(h)+ F_mF_m \subset \Big[0,\frac{m^2+m}{h}\Big]
\cup
\bigcup_{ab\leq h} \big( F_k(ab) + U_{m,1,a}U_{m,1,b} \big).
$$
Write $\tau:=(2^{m-1}-1)2^{2-k-2m}$,
$\sigma:=(2^{m-1}-1)2^{1-k-m}$
and $\lambda:=2^{2-k-2m}-1$.
Then we have 
$0<\tau<\sigma<1$.
Define $C_* :=  (1 + 1/(\sigma-\tau) )/(1-\tau)$.
By \eqref{eq:s3-asymmetric-sum} and Lemma~\ref{lem:s3-single-layer-cover},
$$
\N_\rho\bigg(   \bigcup_{ab\leq h} \big( F_k(ab) + U_{m,1,a}U_{m,1,b} \big)  \bigg)
\leq  2C_m \rho^\lambda \sum_{ ab\leq h : a\geq b} a^{-\tau}b^{-\sigma}
\leq  2C_m C_* \rho^\lambda h^{1-\tau}.
$$
Taking $C:=\max\{m^2+m+1,2C_m C_*\}$, we obtain
$$
\N_\rho\big( F_k (h) +  F_m F_m \big)
\leq \frac{m^2+m+1}{h\rho} +2C_m C_* \rho^\lambda h^{1-\tau}
\leq C \rho^{-1} h^{-1} + C \rho^\lambda h^{1-\tau} .
$$
\end{proof}

We now estimate the covering number of $\T_{2,m}$ using the counting argument in \cite{CJW}.
Recall that $d_{2,m}=2^{3m-1}-2^{m-1}+1$.

\begin{lem}\label{lem:s3-integer-cover}
Let $m\in\Z^+$ with $m\geq2$.
There exists a constant $C>0$, depending only on $m$, such that for every integer $n \geq 2$,
$$
\N_{n^{-d_{2,m}} } (\T_{2,m})\leq Cn^{d_{2,m}-1}.
$$
\end{lem}
\begin{proof}
Fix an integer $n\geq2$.
For each $1 \leq i \leq m$, define 
$$
	G_i := \Big\{ \sum_{j=1}^i \frac{1}{n_j} : n_j \in \Z^+ ,  n_j  \leq n^{2^{j-1}} , 1 \leq j \leq i  \Big\}.
$$
Write $G_0 := \{ 0 \}$.
The set $G_i$ gives an alternative representation of Egyptian fraction-good words in \cite{CJW}.
Then we have $\card G_i \leq n^{2^i-1}$ and
$$
	F_m = \bigcup_{i=0}^m  \big( G_i + F_{m-i}(n^{2^i}) \big).
$$
Write $d:=d_{2,m}$ and $\rho:=n^{-d}$.
Therefore
\begin{equation}\label{eq:s3-finite-cover}
\N_\rho(\T_{2,m})
\leq
\sum_{i=0}^m n^{2^i-1}
\N_\rho\big(F_mF_m+F_{m-i}(n^{2^i})\big).
\end{equation}

Fix $0 \leq i \leq m$ and set $h_i:=n^{2^i}$.
Since $d>2^m$, we have $\rho<1/h_i$.
Applying Lemma~\ref{lem:s3-tail-cover} with $k=m-i$
and using the definition of $d$, we obtain
$$
\begin{aligned}
\N_\rho\big(F_mF_m+F_{m-i}(h_i)\big)
&\leq C_{i,m} \Big( \rho^{-1}{h_i}^{-1} + \rho^{2^{i+2-3m}-1} {h_i}^{1-(2^{m-1}-1)2^{i+2-3m}} \Big)  \\
& = C_{i,m}\Big( n^{d-2^i}+n^{d-2^i-(2^{m-1}-1)(2^i-1)2^{i+2-3m}} \Big)\\
&\leq 2C_{i,m} n^{d-2^i}.
\end{aligned}
$$
Substituting this estimate into \eqref{eq:s3-finite-cover} gives
$$
\N_{n^{-d}}(\T_{2,m})
\leq 2 \sum_{i=0}^m C_{i,m} n^{2^i-1}n^{d-2^i}
= 2 n^{d-1} \sum_{i=0}^m C_{i,m}.
$$
Taking $C:=2\sum_{i=0}^m C_{i,m} $ completes the proof.
\end{proof}

\begin{lem}\label{lem:true-taylor-up-dim}
Let $m\in\Z^+$ with $m\geq2$. Then
$$
\updimB(\T_{2,m})\leq1-\frac{1}{d_{2,m}}.
$$
\end{lem}

\begin{proof}
Write $d:=d_{2,m}$.
For $0<\rho<1$, set $n:=\lceil\rho^{-1/d}\rceil$.
Then $n\geq2$, $n^{-d}\leq\rho$ and $n\leq2\rho^{-1/d}$.
By Lemma~\ref{lem:s3-integer-cover},
$$
\N_\rho(\T_{2,m})
\leq
\N_{n^{-d}}(\T_{2,m})
\leq
C_0n^{d-1}
\leq
2^{d-1}C_0\rho^{1/d-1},
$$
where $C_0>0$ depends only on $m$.
The definition of upper Minkowski dimension gives the desired bound.
\end{proof}

\begin{proof}[Proof of Theorem~\ref{thm:true-taylor}]
For $m=1$, we have $\T_{2,1}=  F +F \cdot F =F_2$.
By \eqref{eq:true-taylor-D}, we have $d_{2,1}=4$. By Theorem~\ref{thm:prod}, we obtain
$$
\dimB(\T_{2,1})=\dimB(F_2) = \frac{3}{4} =1-\frac{1}{d_{2,1}}.
$$

For $m\geq2$, Lemmas~\ref{lem:true-taylor-low-dim} and~\ref{lem:true-taylor-up-dim} give
$$
1-\frac{1}{d_{2,m}}
\leq \lowdimB(\T_{2,m})
\leq \updimB(\T_{2,m})
\leq 1-\frac{1}{d_{2,m}}.
$$
Thus the Minkowski dimension exists and has the asserted value.
\end{proof}

\begin{rem}
The numbers $d_{k,m}$ defined in \eqref{eq:true-taylor-D} satisfy the following recurrence.
For every $m\in\Z^+$,
$
d_{1,m}=2^m,
$
and, for every $k\geq2$,
$$
d_{k,m}
= 2d_{k-1,m} + \big(2^{km}-2^k-1\big)2^{(m-1)k(k-1)/2} +2^{(m-1)(k-1)(k-2)/2}.
$$
\end{rem}

\section{Incomparability of product sets and Taylor sets}

In this section, we compare $F_{m+p}F_q$ and $F_m+F_pF_q$.
In the boundary case $q=1$, the two sets coincide.
More precisely,
$$ F_{m+p} F_q = F_{m+p} = F_m + F_p = F_m + F_p  F_1 . $$
If $p\geq q$, \eqref{eq:dim-equal-prod-taylor} shows that $F_{m+p}F_q$ and $F_m+F_pF_q$ share the same Minkowski dimension.
However, these two sets are incomparable when $p \geq q\geq2$.

\begin{lem}\label{lem:sub-one}
For every $m,p,q \in \Z^+$ with $p,q \geq 2$,
we have $F_{m+p} F_q \not\subset F_m + F_p F_q $.
\end{lem}

\begin{proof}
For every $t \in \Z^+$, we have $F_t \subset [0,t]$ and $\max F_t =t$.
Hence,
$$ 
\max ( F_m + F_p F_q ) = \max F_{m} + \max F_p \cdot \max F_{q} =  m+pq .
$$
On the other hand,
$$ \max ( F_{m+p} F_q ) = \max F_{m+p} \cdot \max F_{q}= (m+p)q> m+pq=\max ( F_m + F_p F_q ) .$$
This implies that $F_{m+p} F_q \not\subset ( F_m + F_p F_q )$.
\end{proof}

To prove that $( F_m + F_p F_q ) \not\subset F_{m+p} F_q $, we first introduce the concept of derived sets.
For any bounded set $A\subset\R$, define its derived set by
$$
\D(A) := \big\{ x\in\R: x \mbox{ is an accumulation point of } A  \big\}.
$$
The following two lemmas are standard results for derived sets.

\begin{lem}\label{lem:acc-1}
Let $A,B\subset\R$ be compact. Then
$$
\D(A+B) = \big( \D(A)+B \big) \cup \big( \D(B)+A \big).
$$
\end{lem}

\begin{proof}
We first show that
$$
\big(\D(A)+B\big)\cup\big( \D(B) +A \big)\subset\D(A+B).
$$
For every $a\in\D(A)$, there exists a sequence of distinct points $a_n\in A$ such that $\lim_{n \to \infty} a_n = a$.
Then for every $b\in B$, $\{a_n+b\}$ is a sequence of distinct points in $A+B$ that converge to $a+b$. 
Thus $\D(A)+B\subset\D(A+B)$.
Similarly, we have $\D(B) +A\subset\D(A+B)$.

Conversely, let $x\in\D(A+B)$. Then there exists a sequence of distinct points
$(a_n+b_n) \in A+B$ such that $a_n\in A$, $b_n\in B$ and $\lim_{n\to\infty}(a_n+b_n)=x$. 
Since $A$ is compact, there exists a subsequence $\{a_{n_i}\}$ such that $\lim_{i\to\infty}a_{n_i}=a$
for some $a\in A$.
Moreover,
$$
	\lim_{i \to\infty} b_{n_i} = \lim_{i \to\infty} (a_{n_i} +b_{n_i}) - \lim_{i \to\infty} a_{n_i} =x-a.
$$
Write $b := x-a$.
We have $b\in B$ since $B$ is compact.

If $a\in\D(A)$ or $b\in\D(B)$, then $x=a+b\in(\D(A)+B)\cup(\D(B)+A)$.
Otherwise, $a\notin\D(A)$ and $b\notin\D(B)$.
Then $a_{n_i}=a$ and $b_{n_i}=b$ for all sufficiently large $i$.
This contradicts the fact that the points $a_n+b_n$ are distinct.
Since $x$ is arbitrary, 
$$
\D(A+B) \subset  (\D(A) +B) \cup  (\D(B) +A).
$$
This proves the lemma.
\end{proof}

\begin{lem}\label{lem:acc-2}
Let $A,B \subset \R$ be compact with $0\in\D(A)\cap\D(B)$. Then
$$
\D(A\cdot B) = \big(\D(A)\cdot B\big)\cup\big(\D(B)\cdot A\big).
$$
\end{lem}

\begin{proof}

We first show that $ 0 \in \D(A\cdot B)$.
Pick $a \in A \setminus \{ 0 \}$ and a sequence of distinct points $b_n\in B$ such that
$\lim_{n\to\infty}b_n=0$.
Then $\{a \cdot b_n\}$ is a sequence of distinct points in $A\cdot B$ converging to $0$.
Therefore $ 0 \in \D(A\cdot B)$.

Now we show that
$$
\big(\D(A)\cdot B\big) \subset \D(A\cdot B).
$$
For every $a\in\D(A)$,
there exists a sequence of distinct points $a_n\in A$ such that $\lim_{n\to\infty}a_n=a$.
For every $b\in B \setminus \{ 0 \}$, 
$\{a_nb\}$ is a sequence of distinct points in $A\cdot B$ converging to $ab$, 
which implies that $ab \in \D(A\cdot B)$.
Together with $0\in\D(A\cdot B)$, this gives $\D(A)\cdot B\subset\D(A\cdot B)$.
Similarly, we obtain $\D(B)\cdot A \subset \D (A \cdot B)$.

Conversely, let $x\in\D(A\cdot B)$. Then there exists a sequence of distinct points
$a_nb_n\in A\cdot B$ such that $a_n\in A$, $b_n\in B$, and $\lim_{n\to\infty}a_nb_n=x$.
Since $A, B$ are compact, there exist convergent subsequences $ \{a_{n_i} \} $ and $\{b_{n_i}\}$ 
with the same index sequence $\{n_i\}$.
Write $a := \lim_{i\to\infty}a_{n_i}$ and $b := \lim_{i\to\infty}b_{n_i}$.
Then we have $a \in A$, $b \in B$ and $x=ab$.
If $a\in\D(A)$, then $x=ab\in\D(A)\cdot B$.
Similarly, $x \in\D(B)\cdot A$ if $b\in\D(B)$.
If $a\notin\D(A)$ and $b\notin\D(B)$, then both $a_{n_i}=a$ and $b_{n_i}=b$ hold for all sufficiently large $i$. 
This contradicts the fact that $ a_n b_n $ are distinct. Therefore,
$$
\D(A\cdot B)\subset
\big(\D(A)\cdot B\big)\cup\big(\D(B)\cdot A\big).
$$
This proves the lemma.
\end{proof}

We now apply the preceding two lemmas to the sets $F_p$.

\begin{lem}\label{lem:derived-form}
For every $p \in \Z^+$,
$
\D(F_p)=F_{p-1}.
$
Moreover, for every $p, q \in \Z^+$,
$$
\D(F_pF_q) = ( F_{p-1} F_q ) \cup ( F_p F_{q-1}).
$$
\end{lem}

\begin{proof}

Since $\D(F)=\{0\}=F_0$, we prove the first assertion by induction. 
For every $p \geq 2$, recall that $F_p=F_{p-1}+F$.
By Lemma~\ref{lem:acc-1}, we obtain
$$
\D(F_p)
= (\D(F_{p-1})+F )\cup (\D(F) + F_{p-1}) 
= ( F_{p-2}+F )\cup (F_0 + F_{p-1} )
= F_{p-1}.
$$

Since $F_p$ and $F_q$ are compact sets and $0\in\D(F_p)\cap\D(F_q),$
Lemma~\ref{lem:acc-2} gives
$$
\D(F_pF_q)
= \big( \D(F_p) \cdot F_q \big) \cup \big(  \D(F_q)  \cdot F_p\big)
= ( F_{p-1} F_q ) \cup ( F_p F_{q-1} ).
$$
This proves the lemma.
\end{proof}

For every integer $t\geq2$, define
$$
\K_t := \bigcup_{a,b\in \Z^+: a+b=t} F_aF_b = \bigcup_{a=1}^{t-1} F_a F_{t-a} .
$$
Then we have $\K_2= F_1 F_1 =F$ and $\K_3=F_2 F_1 \cup F_1 F_2 =F_2$.
Lemma~\ref{lem:derived-form} suggests the following recursive relation for the sets $\K_t$.

\begin{lem}\label{lem:K-derived}
For every $t \geq 3$, we have $\D ( \K_t ) = \K_{t-1}$.
\end{lem}

\begin{proof}
Since the union is finite, the definition of accumulation points gives
$$
\D(\K_t) = \bigcup_{ a=1}^{t-1 } \D(F_aF_{t-a}).
$$
By Lemma~\ref{lem:derived-form}, 
$ \D(F_aF_{t-a}) = F_{a-1}F_{t-a}\cup F_aF_{t-a-1}$ 
for every $1 \leq a \leq t-1$.
Therefore
$$
\D(\K_t) =\bigcup_{ a=1}^{t-1 } \big( F_{a-1}F_{t-a}\cup F_aF_{t-a-1} \big)
= \bigcup_{ a=0}^{t-1 }  F_aF_{t-a-1} 
=F_0 F_{t-1} \cup \K_{t-1}.
$$
Since $F_0=\{0\}$, we have $F_0F_{t-1}=\{0\}\subset\K_{t-1}$.
Therefore $\D(\K_t)=\K_{t-1}$.
\end{proof}

We begin with the smallest nontrivial example $(m,p,q)=(1,2,2)$.

\begin{lem}\label{lem:tmp-k-5}
$(F+F_2F_2) \setminus \K_5$ is non-empty. 
More precisely, $23/ 8 \in (F+F_2F_2) \setminus \K_5$.
\end{lem}

\begin{proof}
First, $23/8 \in (F+F_2F_2) $ follows from 
$$
\frac{23}{8} = 1+\Big(1+\frac{1}{2}\Big)\Big(1+\frac{1}{4}\Big).
$$
Since $\K_5 =F_1 F_4 \cup F_2 F_3 = F_4 \cup F_2 F_3 $, it suffices to check the following two cases.

Assume that $23/8 \in F_4$.
Notice that
$$
	1 + \frac{1}{2} + \frac{1}{2} + \frac{1}{2} =  \frac{5}{2} <  \frac{23}{8} <3.
$$
It follows that exactly two of the four summands are equal to $1$.
Hence $23/8 \in F_4$ would imply that $7/8 \in F_2$.
However, it is easy to check $7/8 \notin F_2$.
This contradicts $23/8 \in F_4$.

Now we assume that $23/8=xy$ with $x \in F_2$ and $y \in F_3$.
If $x<1$, then $x \leq 5/6$ and $y \geq 69/20>3$.
If $x=1$, $ y= 23/8 \notin F_4$ and hence $y \notin F_3$.

If $x>1$, there exists $n \in \Z^+$ such that $x=1+1/n$.
For $n \geq 7$, we have $y \in (5/2,3)$ so that $ y \notin F_3 $.
For $3 \leq n \leq 6$, we have $y \in (2,5/2)$ so that $ y = 2+1/k$ for some $k \in \Z^+$.
A direct check for $n=3,4,5,6$ shows that $y-2$ is not a unit fraction.
For $n=2$, $y=23/12 \in (11/6,2)$.
Since $F_3\cap(11/6,2)=\emptyset$, we have $y\notin F_3$.
For $n=1$, we have $y=23/16$, and it is easy to check that $y\notin F_3$.
In conclusion, $y \notin F_3$ if $x >1$.

As a result, $23/8\notin \K_5$ and the proof is complete.
\end{proof}

Building on the preceding example, we obtain the following result.

\begin{lem}\label{lem:sub-two}

For every $m\geq 1$ and $p,q\geq 2$, there exists
\begin{equation}\label{eq:exist-z}
z \in ( F_m + F_p F_q) \setminus \K_{m+p+q}.
\end{equation}

\end{lem}

\begin{proof}

Lemma~\ref{lem:tmp-k-5} gives the case $(m,p,q)=(1,2,2)$.
We proceed by induction on $m+p+q$.

Assume that \eqref{eq:exist-z} holds for $(m,p,q)$.
Choose $a \in F_m$, $b \in F_p$ and $c \in F_q$ such that
$$
z = a+bc \in ( F_m + F_p F_q ) \setminus \K_{m+p+q}.
$$
Write $t := m+p+q$.
Since $F_m \subset F_{t-1} \subset \K_t $, we must have $b,c>0$.
By Lemma~\ref{lem:K-derived}, we have $\D (\K_{t+1})=\K_t$.
Therefore $z\notin \D(\K_{t+1})$.
By the definition of accumulation points, there exists $\vep>0$ such that
$$
	(z,z+\vep) \cap \K_{t+1} = \emptyset .
$$
We choose $n \in \Z^+$ sufficiently large that $n \vep > \max \{p,q\}$.

For $(m+1,p,q)$, we set
$$
z_1 := a+\frac{1}{n}+bc \in (z,z+\vep).
$$
Since
$a+1 / n \in F_{m+1}$,
we have $z_1 \in (F_{m+1} + F_p F_q)$.
Thus, $ z_1 \notin \K_{t+1}$.

For $(m,p+1,q)$, set
$$
z_2 := a+\Big(b+\frac{1}{n}\Big)c = z+\frac{c}{n}  
\in (z,z+ \vep).
$$
Since $b+1/n \in F_{p+1}$, we have $ z_2 \in ( F_m + F_{p+1} F_{q})$.
Hence $z_2 \notin \K_{t+1}$.
Similarly, we set $z_3 := a+b(c+1/n)$ for $(m,p,q+1)$.
Then we have $z_3 \in ( F_m + F_p F_{q+1})$ and $z_3 \notin \K_{t+1}$.

Since every $(m,p,q)$ with $m\geq1$ and $p,q\geq2$ can be obtained
from $(1,2,2)$ by finitely many such steps, the proof is complete.
\end{proof}

\begin{proof}[Proof of Theorem~\ref{thm:non-equal}]
First,
$F_{m+p} F_q \not\subset F_m + F_p F_q $ follows from Lemma~\ref{lem:sub-one}.
Notice that $F_{m+p} F_q  \subset \K_{m+p+q} $.
By Lemma~\ref{lem:sub-two}, there exists
$$
	z \in \big( (F_m + F_p F_q ) \setminus  \K_{m+p+q} \big) \subset \big( (F_m + F_p F_q)  \setminus  F_{m+p} F_q  \big).
$$
Therefore $( F_m + F_p F_q ) \not\subset F_{m+p} F_q  $.
\end{proof}

The result of Lemma~\ref{lem:derived-form} extends naturally to general product sets $\P_\bfm$.

\begin{rem}\label{rem:derived-mult}
Set $\D^1(A):=\D(A)$.
For every $r \in \Z^+$, define
$$
\D^{r+1}(A):=\D\big(\D^r(A)\big).
$$
Iterating Lemma~\ref{lem:K-derived}, we obtain $\D^r(\K_t)=\K_{t-r}$ for every $1 \leq r\leq t-2$.

Fix $\bfm = (m_1, \ldots,m_s)$. Write $K := \sum_{i=1}^s (m_i-1)$.
For every $r \in \Z^+$, we set
$$
	\Lambda_r := \Big\{ \bfn=(n_1,\ldots,n_s)\in\Sigma^s: \sum_{i=1}^s n_i =K+s-r,  1\leq n_i\leq m_i , 1 \leq i \leq s \Big\}.
$$
If $1 \leq r\leq K$, $\Lambda_r$ is nonempty and 
$$
\D^r(\P_{\bfm}) = \bigcup_{ \bfn \in \Lambda_r } \P_{\bfn}.
$$
In particular, for $r=K$, we have $\Lambda_K = \{ (1,\ldots,1) \}$.
Therefore 
$$
\D^K (\P_{\bfm}) = \P_{(1,\ldots,1)} = F.
$$
Hence 
$\D^{K+1}(\P_{\bfm}) = \D(F) =\{0\}$ and $\D^{K+2}(\P_{\bfm})=\emptyset$.
\end{rem}

\subsection*{Acknowledgements}
\medskip
The authors would like to thank Hai-Peng Chen, Lu-Lu Fang, Jun-Jie Miao, Huo-Jun Ruan, and Yu-Feng Wu for helpful discussions.

\bibliographystyle{amsplain}

\end{document}